# A Note on Banach Principle for *JW*-algebras

GENADY Ya. GRABARNIK [1] AND ALEXANDER A. KATZ [2]

[1] T. J. Watson IBM Research Center
19 Skyline Drive, Hawthorne, NY 10510
USA
genady@us.ibm.com

[2] Department of Mathematics and Computer Science
St. John's University
300 Howard Avenue, DaSilva Hall 314, Staten Island, NY 10301
USA
katza@stjohns.edu

*Abstract:* - In the sequel we establish the Banach Principle for semifinite *JW*-algebras without direct summand of type $I_2$, which extends the recent results of Chilin and Litvinov on the Banach Principle for semifinite von Neumann algebras to the case of *JW*-algebras.

*Key-Words:* - von Neumann algebras, Jordan operator algebras, *JW*-algebras, Banach Principle, *-algebra of $\tau$-measurable operators affiliated to a semifinite von Neumann algebra, Jordan algebra of $\tau$-measurable operators affiliated to a semifinite *JW*-algebra.

## 1 Introduction

Let $(\Omega, \Sigma, \mu)$ be a probability space. Denote by $\mathcal{L} = \mathcal{L}(\Omega, \mu)$ the set of all (classes of) complex-valued measurable functions on $\Omega$. Let $\tau_\mu$ be the measure topology on $\mathcal{L}$. The classical Banach Principle (see for example [13]) can be stated as follows:

**Classical Banach Principle.** *Let $(X \parallel \cdot \parallel)$ be a Banach space, and let $a_n : (X \parallel \cdot \parallel) \to (\mathcal{L}, \tau_\mu)$ be a sequence of continuous linear maps. Consider the following properties:*

*(I) the sequence $\{a_n(x)\}$ converges almost everywhere (a.e.) for every $x \in X$;*

*(II) $\hat{a}(x)(\omega) = \sup_n |a_n(x)(\omega)| < \infty$ a.e. for every $x \in X$;*

*(III) (II) holds, and the maximal operator $\hat{a} : (X \parallel \cdot \parallel) \to (\mathcal{L}, \tau_\mu)$ is continuous at $0$;*

*(IV) the set $\{x \in X : \{a_n(x)\} \text{ converges a.e.}\}$ is closed in $X$.*

*Then the implications $(I) \Rightarrow (II) \Rightarrow (III) \Rightarrow (IV)$ are always true. If in addition, there exists a dense subset $D \subset X$, such that the sequence $\{a_n(x)\}$ converges a.e. for every $x \in D$, then all four conditions $(I) - (IV)$ above are equivalent.*

The Banach Principle above was often applied in the case $X = (L^p, \parallel \cdot \parallel_p)$, where $1 \le p < \infty$. However, in the case $p = \infty$ the uniform topology on $L^\infty$ appears to be too strong for the classical Banach Principle to be effective in $L^\infty$. For example, one can notice that continuous functions are not uniformly dense in $L^\infty$.

Bellow and Jones [9], using the fact that the unit ball $L_1^\infty = \{x \in L^\infty : \parallel x \parallel_\infty \le 1\}$ is complete in $\tau_\mu$, suggested to consider the measure topology on $L^\infty$ by replacing $(X \parallel \cdot \parallel)$ by $(L_1^\infty, \tau_\mu)$. Since $L_1^\infty$ is not a linear space, geometrical complications occur, which, however, were successfully resolved in [9].

Non-commutative versions of Banach Principle for measurable operators affiliated to a semifinite von Neumann algebra were established in [14] and [18]. These results were extended to the case of semifinite *JBW*-algebras in [16] and [17], following ideas introduced in [1], [2], [4], [5], [10]. A non-commutative version of the Banach Principle for $L^\infty$ was proposed by Chilin and Litvinov in [12]. The present notes are devoted to a presentation of an extension of the results in [12] to the case of JW-algebras without direct summand of the type $I_2$.

## 2 Preliminaries

Let $M$ be a semifinite von Neumann algebra of bounded operators acting on a complex Hilbert space $H$ ([11]), and let $B(H)$ be the algebra of all bounded operators on $H$. A densely defined closed operator $x$ on $H$ is called *affiliated* to $M$ if $y'z = zy'$, with $z \in M$ ([19], [21]). Denote by $P(M)$ the complete lattice of projections in $M$. Let $\tau$ be a faithful normal semifinite trace on $M$. Denote by $e^\perp = \mathbf{1} - e$ the orthogonal complemented projection for the projection $e \in P(M)$. An operator $x$ affiliated to M is called $\tau$-*measurable* if $\forall \varepsilon > 0$, $\exists e \in P(M)$ with $\tau(e^\perp) \leq \varepsilon$ such that $eH$ belongs to the domain of the operator $x$. Let $L(M, \tau)$ be the set of all $\tau$-measurable operators affiliated to $M$.

Set $V(\varepsilon, \delta) = \{x \in L(M, \tau) : \| xe \| < \delta$ for some $e \in P(M)$ with $\tau(e^\perp) < \varepsilon\}$, for arbitrary $\varepsilon > 0$ and $\delta > 0$, where $\|\cdot\|$ stands for the operator norm on $B(H)$. The topology $t_\tau$ defined on $L(M, \tau)$ by the family $\{V(\varepsilon, \delta) : \varepsilon > 0, \delta > 0\}$ of neighborhoods of zero is called the *measure topology* ([19], [21]).

**Theorem 1.** $(L(M, \tau), t_\tau)$ *is a complete metrizable topological *-algebra.*

**Proof.** *See [19], [21] for details.* □

**Proposition 1.** *For any $d > 0$, the sets*
$M_d = \{x \in M : \| x \| \leq d\}$, $M_d^h = \{x \in M_d : x = x^*\}$
*are $t_\tau$-complete.*

**Proof.** *See [12] for details.* □

A sequence $\{y_n\} \subset L(M, \tau)$ is said to converge *almost uniformly* (a.u.) to $y \in L(M, \tau)$ if $\forall \varepsilon > 0$, $\exists e \in P(M)$ with $\tau(e^\perp) < \varepsilon$ such that $\|(y - y_n)e\| \to 0$.

**Proposition 2.** *For $\{y_n\} \subset L(M, \tau)$ the conditions*

*(i) $\{y_n\}$ converges a.u. in $L(M, \tau)$;*

*(ii) $\forall \varepsilon > 0$, $\exists e \in P(M)$ with $\tau(e^\perp) < \varepsilon$ such that $\|(y_m - y_n)e\| \to 0$ as $m, n \to \infty$;*

*are equivalent.*

**Proof.** *See [12] for details.* □

The following theorem is a non-commutative version of Riesz theorem ([13]).

**Theorem 2.** *If $\{y_n\} \subset L(M, \tau)$ and $y = t_\tau - \lim_{n \to \infty} y_n$, then $y = a.u. - \lim_{k \to \infty} y_{n_k}$ for some subsequence $\{y_{n_k}\} \subset \{y_n\}$.*

**Proof.** *See [21] and [14] for details.* □

Let $A$ be a semifinite JW-subalgebra of $B(H)_{SA}$ without a direct summand of type $I_2$ (see [15] and [20] for definitions), $P(A)$ be the complete lattice of projections in $A$, and $\tau$ be a faithful normal semifinite trace on A. Let $M = M(A)$ be the von Neumann enveloping algebra of the Jordan algebra $A$. Then $\tau$ can be uniquely extended to a faithful normal semifinite trace on $M$, for which we will use the same symbol $\tau$ (see [3], [6] and [8]). A self adjoint operator $x \in L(M, \tau)$ is called affiliated to a JW-algebra $A$, if all its spectral projections belong to $A$ (see [3], [6], [7] and [8]). An operator $x$ affiliated to A is called $\tau$-*measurable* if $\forall \varepsilon > 0$, $\exists e \in P(A)$ with $\tau(e^\perp) \leq \varepsilon$ such that $eH$ belongs to the domain of the operator $x$. Let $L(A, \tau)$ be the set of all $\tau$-measurable operators affiliated to $A$.

**Proposition 3.** *An operator $x \in L(M, \tau)_{SA}$ is affiliated to A iff $x \in L(A, \tau)$.*

**Proof.** *Follows from arguments in [20].* □

**Theorem 3.** $(L(A, \tau), t_\tau)$ *is a complete topological Jordan subalgebra of $(L(M, \tau), t_\tau)_{SA}$.*

**Proof.** *A direct consequence of Theorem 1 and arguments in [20].* □

A sequence $\{y_n\} \subset L(A, \tau)$ is said to converge *bilaterally with square almost uniformly* (b.s.a.u.) to $y \in L(A, \tau)$ if $\forall \varepsilon > 0$, $\exists e \in P(A)$ with $\tau(e^\perp) < \varepsilon$ such that $\|e(y - y_n)^2 e\| \to 0$.

**Proposition 4.** *For $\{y_n\} \subset L(A, \tau) \subset L(M, \tau)_{SA}$ the conditions:*

*(i) $\{y_n\}$ converges a.u. in $L(M, \tau)$;*

*(ii) $\forall \varepsilon > 0$, $\exists e \in P(M)$ with $\tau(e^\perp) < \varepsilon$ such that $\|(y_m - y_n)e\| \to 0$ as $m, n \to \infty$;*

*(iii) $\{y_n\}$ converges b.s.a.u. in $L(A, \tau)$;*

*(iv) $\forall \varepsilon > 0$, $\exists e \in P(A)$ with $\tau(e^\perp) < \varepsilon$ such that $\|e(y_m - y_n)^2 e\| \to 0$ as $m, n \to \infty$;*

*are equivalent.*

**Proof.** *From* $\|e(y_m - y_n)^2 e\| =$
$= \|e(y_m - y_n)(y_m - y_n)e\| =$

$$=\|((y_m - y_n)e)^*(y_m - y_n)e\| \leq$$
$$\leq \|((y_m - y_n)e)^*\| \cdot \|(y_m - y_n)e\| =$$
$$= \|(y_m - y_n)e\|^2, \text{ so we can see that}$$

b.s.a.u. fundamentalness of a sequence in a reversible JW-algebra ([15], [8]) is equivalent to a.u. fundamentalness of the same sequence in its von Neumann enveloping algebra $M = M(A)$. Thus the statement follows from Proposition 2 above. □

The Riesz theorem 2 above will take the following form.

**Theorem 4.** *If* $\{y_n\} \subset L(A, \tau)$ *and*
$y = t_\tau - \lim_{n \to \infty} y_n$, *then* $y = b.s.a.u. - \lim_{k \to \infty} y_{n_k}$
*for some subsequence* $\{y_{n_k}\} \subset \{y_n\}$.

**Proof.** *Directly follows from Proposition 4 and Theorem 2 above.* □

## 3 Bilateral with square uniform equicontinuity for sequences of maps into $L(A,\tau)$

Let E be an arbitrary set. If $a_n : E \to L(A, \tau)$, $x \in E$, and $b \in A$ such that $\{b(a_n(x))^2 b\} \subset A$. Denote $S(\{a_n^2\}, x, b) = \sup_n \|b(a_n(x))^2 b\|$.

The following Lemma is valid.

**Lemma 1.** *Let* $(X, +)$ *be a semigroup, and* $a_n : X \to L(A, \tau)$ *be a sequence of additive maps. Assume that* $\bar{x} \in X$ *is such that* $\forall \varepsilon > 0$, $\exists \{x_k\} \subset X$, *and* $p \in P(A)$ *with* $\tau(p^\perp) < \varepsilon$, *such that :*
*(i)* $\{a_n(\bar{x} + x_k)\}$ *converges b.s.a.u. as* $n \to \infty$, *for every* $k \in N$;
*(ii)* $S(\{a_n^2\}, x_k, p) \to 0$, *as* $k \to \infty$.
*Then the sequence* $\{a_n(\bar{x})\}$ *converges b.s.a.u. in* $L(A, \tau)$.

**Proof.** *Follows from [12] and Proposition 4.* □

Let $(X, t)$ be a topological space, and $a_n : X \to L(A, \tau)$ and $x_0 \in X$ be such that $a_n(x_0) = 0$ for $n \in N$. A sequence $\{a_n\}$ is called bilaterally with square equicontinuous at $x_0$ if $\forall \varepsilon, \delta > 0$, $\exists$ a neighborhood $U$ of $x_0$ in $(X, t)$ such that $a_n U \subset V(\varepsilon, \delta) \cap L(A, \tau)$, $n \in N$, i.e. $\forall x \in U$ and $\forall n \in N$ one can find a projection $e = e(x, n) \in P(A)$ with $\tau(e^\perp) < \varepsilon$, satisfying $\|e(a_n(x))^2 e\| < \delta$.

Let now $x_0 \in E \subset X$. A sequence $\{a_n\}$ is called **bilaterally with square uniformly equicontinuous** at $x_0$ on E, if $\forall \varepsilon, \delta > 0$, $\exists$ a neighborhood $U$ of $x_0$ in $(X, t)$ such that $\forall x \in E \cap U$, $\exists e = e(x) \in P(A)$ with $\tau(e^\perp) < \varepsilon$, satisfying $S(\{a_n^2\}, x, e) < \delta$.

**Proposition 5.** *Let the sequence* $\{a_n\}$ *and* $x_0 \in E \subset X$ *be as above. Then,*
*(i)* $\{a_n\}$ *is equicontinuous at* $x_0$ *on E into* $L(M, \tau)$ *iff it is bilaterally with square equicontinuous at* $x_0$ *on E into* $L(A, \tau)$;
*(ii)* $\{a_n\}$ *is uniformly equicontinuous ([12]) at* $x_0$ *on E into* $L(M, \tau)$ *iff it is bilaterally with square uniformly equicontinuous at* $x_0$ *on E into* $L(A, \tau)$.

**Proof.** *Directly follows from Proposition 4 and arguments in [12].* □

Theorem 1 and theorem 3 established that that $(L(M, \tau), t_\tau)$ is a complete metrizable topological *-algebra, and $(L(A, \tau), t_\tau)$ is a complete metrizable topological Jordan subalgebra of $(L(M, \tau), t_\tau)_{SA}$. In [12] it has been established that for any $d > 0$, the sets
$$M_d = \{x \in M : \|x\| \leq d\}, \text{ and}$$
$$M_d^h = \{x \in M_d : x = x^*\} \text{ are } t_\tau \text{-comp-}$$
lete. It is easy to see that the set $A_d = M_d^h \cap A$ is $t_\tau$-complete too.

**Lemma 2.** *Let* $d > 0$. *If* $a_n : A \to L(A, \tau)$ *be a sequence of additive maps. Then it is bilaterally with square uniformly equicontinuous at* 0 *on* $A_d$ *iff it is uniformly equicontinuous at* 0 *on* $M_d$ *(where in the second condition we mean that all maps are extended by linearity to the sequence of additive maps* $M \to L(M, \tau)$ *).*

**Proof.** *Directly follows from Proposition 5, and arguments in [12] and [8].* □

**Lemma 3.** *Let a sequence* $a_n : A \to L(A, \tau)$ *of additive maps be bilaterally with square uniformly equicontinuous at* 0 *on* $A_d$ *for some* $0 < d \in R$. *Then* $\{a_n\}$ *is as well bilaterally with square* uniformly equicontinuous at $0$ on $A_s$ for every $0 < s \in \mathbf{R}$.
**Proof.** *Directly follows from Lemma 2 and arguments in [12] and [8].* □

## 4 Main results

Let $0 \in E \subset A$. For a sequence $a_n : (A, t_\tau) \to L(A, \tau)$, consider the following conditions:
- Bilateral with square almost uniform convergence of $\{a_n(x)\}$ for every $x \in E$ (BSCNV (*E*));
- Bilateral with square uniform equicontinuity at 0 on *E* (BSCNT (*E*));
- Closedness in $(E, t_\tau)$ of the set $C(E) = \{x \in E : \{a_n(x)\}$ converges b.s.a.u.$\}$ (BSCLS (*E*)).

In this section we will discuss relationships among the conditions (BSCNV ($A_1$)), (BSCNT ($A_1$)), and (BSCLS ($A_1$)).

**Theorem 5.** *Let $a_n : A \to L(A, \tau)$ be a (BSCNV ($A_1$)) sequence of positive $t_\tau$-continuous linear maps with $a_n(1) \leq 1$, $n \in N$. Then the sequence $\{a_n\}$ is also (BSCNT ($A_1$)).*
**Proof.** *Directly follows from arguments in [12] and the previous section.* □

**Theorem 6.** *A (BSCNT ($A_1$)) sequence of additive maps $a_n : A \to L(A, \tau)$ is as well (BSCLS ($A_1$)).*
**Proof.** *Directly follows from arguments in [12] and the results of the previous section.* □

**Theorem 7.** *Let $a_n : A \to L(A, \tau)$ be a sequence of positive $t_\tau$-continuous linear maps such that $a_n(1) \leq 1$, $n \in N$. If a sequence $\{a_n\}$ is (BSCNV (D)) with D being $t_\tau$-dense in $A_1$, the conditions (BSCNV ($A_1$)), (BSCNV ($A_1$)), and (BSCLS ($A_1$)) are equivalent.*
**Proof.** *Directly follows from [12] and the results of the previous section.* □

## 5 Conclusion

Results of the present notes extend the results of [12] to the case of *JW*-algebras without direct summand of type $I_2$. In a new manuscript under preparation we extend these results to the case of bilateral almost uniform convergence on semifinite von Neumann algebras and semifinite *JBW*-algebras without direct summand of type $I_2$.